\documentclass[12pt]{article} 
\usepackage{amssymb,latexsym,theorem,hyperref}
\makeatletter
\newcommand{\Ab}{\mathop{\operator@font{Abutment}}\nolimits}
\newcommand{\Diag}{\mathop{\operator@font{Diag}}\nolimits}
\newcommand{\diag}{\mathop{\operator@font{diag}}\nolimits}
\newcommand{\Proj}{\mathop{\operator@font{Proj}}\nolimits}
\newcommand{\Spec}{\mathop{\operator@font{Spec}}\nolimits}
\newcommand{\res}{\mathop{\operator@font{res}}\nolimits}
\newcommand{\red}{{\mathop{\operator@font{red}}\nolimits}}
\newcommand{\ind}{\mathop{\operator@font{ind}}\nolimits}
\newcommand{\hull}{\mathop{\operator@font{hull}}\nolimits}
\newcommand{\rank}{\mathop{\operator@font{rank}}\nolimits}
\newcommand{\even}{{\operator@font{even}}}
\newcommand{\hgt}{\mathop{\operator@font{ht}}\nolimits}
\newcommand{\gr}{\mathop{\operator@font{gr}}\nolimits}
\newcommand{\im}{\mathop{\operator@font{im}}\nolimits}
\newcommand{\st}{\mathop{\operator@font{St}}\nolimits}
\newcommand{\GL}{\mathop{\operator@font{GL}}\nolimits}
\newcommand{\Gr}{\mathop{\operator@font{Gr}}\nolimits}
\newcommand{\gl}{\mathop{\mathfrak{gl}}\nolimits}

\newcommand{\SL}{\mathop{\operator@font{SL}}\nolimits}
\newcommand{\Ext}{\mathop{\operator@font{Ext}}\nolimits}
\newcommand{\Hom}{\mathop{\operator@font{Hom}}\nolimits}
\newcommand{\Pic}{\mathop{\operator@font{Pic}}\nolimits}
\newcommand{\Z}{\mathbb Z}
\newcommand{\Q}{\mathbb Q}
\newcommand{\R}{\mathbb R}

\newcommand{\Gm}{{\mathbb G_m}}
\newcommand{\cE}{\mathcal E}
\newcommand{\cF}{\mathcal F}
\newcommand{\cG}{\mathcal G}
\newcommand{\cO}{\mathcal O}
\newcommand{\cS}{\mathcal S}
\newcommand{\cQ}{\mathcal Q}
\newcommand{\cHom}{\mathcal Hom}
\newcommand{\m}{\mathfrak m}

\newcommand{\cM}{\mathcal M}
\newcommand{\N}{\mathcal N}

\newcommand{\AG}{{$AG$}}
\newcommand{\CG}{{$CG$}}
\newcommand{\GD}{{$G\Delta$}}
\newcommand{\AGD}{{$AG\Delta$}}
\newcommand{\CGD}{{$CG\Delta$}}
\newcommand{\CGL}{{$CG\Lambda$}}
\newcommand{\qed}{\unskip\nobreak\hfill\hbox{ $\Box$}}
\makeatother

\newtheorem{Proposition}[subsection]{Proposition}
\newtheorem{Theorem}[subsection]{Theorem}
\newtheorem{Lemma}[subsection]{Lemma}
\newtheorem{Corollary}[subsection]{Corollary}

\theorembodyfont{\normalfont}
\newtheorem{Remark}[subsection]{Remark}

\newtheorem{Example}[subsection]{Example}
\newtheorem{Definition}[subsection]{Definition}
\newtheorem{Notation}[subsection]{Notation}
\newtheorem{Assumption}[subsection]{Assumption}

\begin{document}
\title{Finite Schur filtration dimension for modules
over an algebra with Schur filtration}
\author{ Vasudevan Srinivas and Wilberd van der Kallen}
\date{}
\maketitle
\sloppy
\begin{abstract}Let $G=\GL_N$ or $\SL_N$ as
reductive linear algebraic group over a field $k$ of
characteristic $p>0$. We prove several results that were previously
established only when $N\leq5$ or $p>2^N$:
Let $G$ act rationally on a finitely
generated commutative $k$-algebra $A$ and let $\gr A$ be the Grosshans
graded ring. 
We show that the cohomology algebra $H^*(G,\gr A)$ is finitely
generated over $k$. If moreover $A$ has a good filtration and
$M$ is a noetherian $A$-module with compatible $G$ action,
then $M$ has finite good filtration dimension and the
$H^i(G,M)$ are noetherian $A^G$-modules. To obtain results in this
generality, we employ functorial resolution of the ideal of the diagonal in
a product of Grassmannians. 
\end{abstract}

\section{Introduction}
Consider  a connected reductive linear algebraic  group
$G$ defined over a field $k$ of positive characteristic $p$.
We say that $G$ has the cohomological finite generation property (CFG)
if the following holds:
Let $A$ be a finitely generated commutative $k$-algebra on which $G$ acts
 rationally by $k$-algebra automorphisms. (So $G$ acts from the right
 on $\Spec(A)$.)
Then the cohomology ring $H^*(G,A)$ is finitely generated
as a $k$-algebra. Here, as in
\cite[I.4]{Jantzen}, we use the cohomology introduced by
Hochschild, also known as `rational cohomology'.

The intent of this paper is to take one more step towards proving the 
 conjecture that every reductive
linear algebraic group has property (CFG). 
The proof will be finished by Antoine Touz\'e, cf.\ \cite{Touze}. 
The key point of the present work is to remove restrictions on the characteristic
from \cite{fgfdim}.  

Our proofs use resolution of the diagonal in products of Grassmannians.
Thus they apply only to the groups $\SL_N$, $\GL_N$.
But  recall (\cite{cohGrosshans}, \cite{reductive}, \cite{fgfdim}) that
for the conjecture these cases suffice. Also recall that the conjecture 
implies the main results of this paper,
as well as their analogues for other reductive groups.

To formulate the main results, let $N\geq1$ and let $G$ be the connected
reductive linear algebraic group 
$\GL_N$ or $\SL_N$ over an algebraically closed field $k$ of
characteristic $p>0$.
Let $A$ be a finitely generated commutative $k$-algebra on which $G$ acts
rationally by $k$-algebra automorphisms. Let $M$ be a noetherian $A$-module
on which $G$ acts compatibly. This means that the structure map $A\otimes M\to
M$ is a $G$-module map.
Our main theorem is

\begin{Theorem}\label{maingood}
If $A$ has a good filtration, then $M$ has finite good filtration dimension and
each $H^i(G,M)$ is a noetherian $A^G$-module.
\end{Theorem}

One may also formulate the first part in terms of polynomial representations of $\GL_N$.
Recall that a finite dimensional (as $k$ vector space)
rational representation of $\GL_N$ is called polynomial if
it extends to the monoid of $N$ by $N$ matrices without poles along the 
locus where the determinant vanishes.
 Unlike Green \cite {Green} we
cannot restrict ourselves to finite dimensional representations, so we define
a representation to be polynomial if it is a union of finite dimensional
polynomial representations. In other words, we allow infinite dimensional comodules
for the 
bi
algebra of regular functions on the monoid. 

So let $A$ be a finitely generated commutative $k$-algebra on which $\GL_N$ acts
polynomially by $k$-algebra automorphisms. Let $M$ be a noetherian $A$-module
on which $\GL_N$ acts compatibly and polynomially.

\begin{Theorem}\label{mainSchur}
If $A$ has Schur filtration, then $M$ has finite Schur filtration dimension.
\end{Theorem}

\begin{Remark}
The $H^i(\GL_N,M)$ are less interesting now, because the part of nonzero
polynomial degree in $M$ does not contribute to $H^i(\GL_N,M)$.
\end{Remark}

Now let  $A$ 
be a finitely generated commutative $k$-algebra on which $\SL_N$ acts
rationally by $k$-algebra automorphisms. One then has a Grosshans
graded algebra $\gr A$ and we can remove the restrictions on the characteristic
in \cite[Theorem 1.1]{cohGrosshans}:
\begin{Corollary}\label{mainGrosshans}
\label{fingrgr}
The $k$-algebra $H^*(\SL_N,\gr A)$ is finitely generated.
\end{Corollary}

The method of proof of the main result
is based on the functorial resolution \cite{twisted}
of the diagonal of
$Z\times Z$  when $Z$ is a Grassmannian of subspaces of $k^N$.
This is used inductively to study equivariant sheaves on a product $X$
of such Grassmannians.
That leads to a special case of the theorems, with $A$ equal to the 
Cox ring of $X$, multigraded by the Picard group $\Pic(X)$, and $M$ compatibly
 multigraded.
Next one treats cases when on the same $A$ the multigrading is replaced 
with a `collapsed' grading with smaller value group
and $M$ is only required to be 
 multigraded compatibly with this new grading. Here the trick is that 
an associated graded of $M$ has a multigrading that is collapsed a little less.
The suitably multigraded Cox rings now replace the 
`graded polynomial algebras with good filtration' of \cite{cohGrosshans}
and the method of \cite{fgfdim} applies to finish the proof of Theorem 
\ref{maingood}. Then Corollary \ref{mainGrosshans} follows in the manner
of \cite{cohGrosshans}.

\paragraph{Acknowledgements}
Our collaboration got started thanks to the 60th birthday conferences for
V.~B.~Mehta and S.~M.~Bhatwadekar at TIFR Mumbai in 2006. 
Much of the subsequent work was done at the university of Bielefeld, which we thank for
its hospitality.

\section{Recollections and conventions}
Some unexplained notations, terminology, properties, \ldots can be
found in \cite{Jantzen}.

\label{G=GL}
{}From now on, with the exception of section \ref{G=SL}, we put
 $G=\GL_N$, with $B^+$ its subgroup of upper triangular matrices,
$B^-$ the subgroup of lower triangular matrices, $T=B^+\cap B^-$
the diagonal subgroup, $U=U^+$ the unipotent radical of $B^+$.
The roots of $U$ are positive.
The character group $X(T)$ has a basis $\epsilon_1$ \ldots, $\epsilon_N$
with $\epsilon_i(\diag(t_1,\ldots,t_N))=t_i$. 
An element $\lambda=\sum_i \lambda_i\epsilon_i$ of $X(T)$ is often denoted 
$(\lambda_1,\ldots,\lambda_N)$. It is called a polynomial weight if the
$\lambda_i$ are nonnegative. It is called a dominant weight if $\lambda_1\geq
\cdots\geq\lambda_N$. It is called anti-dominant if $\lambda_1\leq
\cdots\leq\lambda_N$. The fundamental weights $\varpi_1$, \dots , $\varpi_N$
are given by $\varpi_i=\sum_{j=1}^i\epsilon_j$.
If $\lambda\in X(T)$ is dominant, then $\ind_{B^-}^G(\lambda)$ is the
\emph{dual Weyl module} or \emph{costandard module} $\nabla_G(\lambda)$,
or simply $\nabla(\lambda)$, with highest weight $\lambda$. 
The \emph{Grosshans height} of $\lambda$ is
$\hgt(\lambda)=\sum_i(N-2i+1)\lambda_i$. It extends to a homomorphism
$\hgt:X(T)\otimes\Q\to\Q$.
The determinant representation has weight $\varpi_N$ and one has
$\hgt(\varpi_N)=0$.
Each positive root $\beta$ has $\hgt(\beta)>0$.
If $\lambda$ is a dominant polynomial weight, then $\nabla_G(\lambda)$
is called a \emph{Schur module}. If $\alpha$ is a partition with at most $N$ parts
then we may view it as a dominant polynomial weight and
the \emph{Schur functor} $S^\alpha$  maps $\nabla_G(\varpi_1)$ to 
$\nabla_G(\alpha)$. (This is the convention followed in \cite{twisted}.
In \cite{Akin} the same Schur functor is labeled with the conjugate
partition $\tilde \alpha$. See also \cite[Thm.~(4.8f), 5.6]{Green}.)
The formula
$\nabla(\lambda)=\ind_{B^-}^G(\lambda)$ just means that $\nabla(\lambda)$
is obtained from the Borel-Weil construction:
$\nabla(\lambda)$ equals $H^0(G/B^-,{\mathcal L}_\lambda)$ for a certain
 line bundle ${\mathcal L}_\lambda$ on the
flag variety $G/B^-$.
There are similar conventions for $\SL_N$-modules. For instance, the costandard
modules for $\SL_N$ are the restrictions of those for $\GL_N$.
The Grosshans height on $X(T)$ induces one on $X(T\cap\SL_N)\otimes\Q$. 
The \emph {multicone} $k[\SL_N/U]$ consists of the $f$ in the coordinate ring
$k[\SL_N]$ that satisfy $f(xu)=f(x)$ for $u\in U\cap \SL_N$. As an 
$\SL_N$-module it is the direct sum of all costandard modules. It is also
a finitely generated algebra \cite{Kempf Ramanathan}, \cite{Grosshans contr}.

\begin{Definition}
A \emph {good filtration} of a $G$-module $V$ is a filtration
$0=V_{\leq -1}\subseteq V_{\leq 0} \subseteq V_{\leq 1}\ldots$ by
$G$-submodules $V_{\leq i}$ with $V=\cup_iV_{\leq i}$, so
that its associated graded $\gr V$ is a direct sum of costandard modules.
A Schur filtration of a polynomial $\GL_N$-module $V$ is a filtration
$0=V_{\leq -1}\subseteq V_{\leq 0} \subseteq V_{\leq 1}\ldots$ by
$\GL_N$-submodules with $V=\cup_iV_{\leq i}$, so
that its associated graded $\gr V$ is a direct sum of Schur modules.
The \emph{Grosshans filtration} of $V$ is the filtration with
$V_{\leq i }$ the largest $G$-submodule of $V$ whose weights
$\lambda$ all satisfy $\hgt(\lambda)\leq i$. Good filtrations and
Grosshans filtrations for $\SL_N$-modules are defined similarly. 
The literature contains more restrictive definitions 
of good/Schur filtrations. Ours are the right ones
when dealing with infinite dimensional representations 
\cite{vdkallen book},
cf.\ \cite[II.4.16 Remark 1]{Jantzen}.
\end{Definition}

\begin{Proposition}\label{Schurgood}
Let $V$ be a polynomial representation of $\GL_N$.
The following are equivalent
\begin{enumerate}
\item $V$ has a good filtration,
\item  $V$ has a Schur filtration,
\item The Grosshans filtration of\/ $V$ is a Schur filtration,
\item The restriction $\res^{\GL_N}_{\SL_N}V$ has a good filtration,
\item  The Grosshans filtration of the restriction $\res^{\GL_N}_{\SL_N}V$ 
is a good filtration,
\item $H^1(\SL_N,k[\SL_N/U]\otimes V)=0$.
\end{enumerate}
\end{Proposition}

\paragraph{Proof}
$3\Rightarrow2\Rightarrow1\Rightarrow4\Rightarrow6\Rightarrow5$
is well known \cite[II 4.16, proof of A.5]{Jantzen}, 
compare \cite[Exercise 4.1.3]{vdkallen book}.
Now assume 5. We may decompose $V$ into weight spaces
(also known as polynomial degrees) for the center of $G$. One may
replace $V$ by one of these weight spaces. The Grosshans filtration
of   $\res^{\GL_N}_{\SL_N}V$ is then a good filtration which may
be reinterpreted as
a Schur filtration on $V$.
\qed

\begin{Definition}
If $V$ is a $\GL_N$-module, and $m\geq-1$ is an integer so that
$H^{m+1}(\SL_N,k[\SL_N/U]\otimes \res^{\GL_N}_{\SL_N} V)=0$, then we say
that $V$ has \emph{good filtration dimension} at most $m$. 
(Compare \cite{Friedlander-Parshall}.)
The case $m=0$ corresponds with $V$ having a good filtration.
And for $m\geq0$ it means that $V$ has a resolution
$$0\to V\to N_0 \to \cdots \to N_m\to 0$$ in which the $N_i$ have good
filtration.
We say that $V$ has good filtration dimension precisely $m$,
notation $\dim_\nabla(V)=m$,
 if $m$ is
minimal so that $V$ has good filtration dimension at most $m$.
In that case $H^{i+1}(\SL_N,k[\SL_N/U]\otimes \res^{\GL_N}_{\SL_N} V)=0$
for all  $i\geq m$. In particular $H^{i+1}(G,V)=0$ for $i\geq m$.
If there is no finite $m$ so that  $\dim_\nabla(V)=m$, then we put
$\dim_\nabla(V)=\infty$. Similar definitions apply to $\SL_N$-modules.

If $V$ is a polynomial representation then $\dim_\nabla(V)$ is also
called the Schur filtration dimension. Indeed if for such $V$ one has
$ \dim_\nabla(V)\leq m$,
$m\geq0$,
then $V$ has  a resolution
$$0\to V\to N_0 \to \cdots \to N_m\to 0$$ in which the $N_i$ have Schur
filtration. 
\end{Definition}

\section{Gradings}\label{gradings}
Let $\Delta=\Z^r$ with standard basis $e_1$, \ldots , $e_r$.
We partially order $\Delta$ by declaring that $I\geq J$ if $I_q\geq J_q$
for $1\leq q\leq r$. The \emph{diagonal} $\diag(\Delta)$ consists of the integer
multiples of the vector $E=(1,\ldots,1)$.
By a \emph{good $G$-algebra} we mean a 
finitely generated commutative $k$-algebra $A$ on which $G$ acts rationally
by $k$-algebra automorphisms so that $A$ has a good filtration as a $G$-module.
We say that $A$ is a \emph{good \GD-algebra} if moreover $A$ is $\Delta$-graded
by $G$-submodules,
$$A=\bigoplus_{I\in\Delta,~I\geq0}A_I$$
 with 
\begin{itemize}
\item
 $A_IA_J\subset A_{I+J}$,
\item $A$
is generated over $A_0$ by the $A_{e_q}$,
\item $G$ acts trivially on $A_0$.
\end{itemize}
Motivated by the Segre embedding we define 
$$\diag(A)=\bigoplus_{I\in\diag(\Delta)}A_I$$ and $\Proj(A):=\Proj(\diag(A))$.
By an \AG-module we will mean a noetherian $A$-module $M$
with compatible $G$-action.
If moreover $M$ is $\Delta$-graded by $G$-submodules $M_I$ so that 
$A_IM_J\subset M_{I+J}$, then we call $M$ an \emph{\AGD-module}.

\begin{Definition}
We call an \AG-module $M$ \emph{negligible} if $M$ has finite 
good filtration dimension
and each $H^i({\SL_N},M)$ is a noetherian $A^{\SL_N}$-module. 
Let $\N$ be the class of the negligible \AG-modules.
\end{Definition}

\begin{Lemma}\label{two/three}
$\N$ has the two out of three property:
If $$0\to M'\to M\to M''\to 0$$ is exact,
and two of $M'$, $M$, $M''$
are negligible, then so is the third.
\end{Lemma}

\paragraph{Proof}
The short exact sequence of Hochschild complexes 
\cite[I.4.14]{Jantzen}
$$0\to C^*(\SL_N,M')\to C^*(\SL_N,M)\to C^*(\SL_N,M'')\to 0$$
is a bicomplex of $A^{\SL_N}$-modules, so the long exact sequence 
$$\cdots\to H^i({\SL_N},M')\to H^i({\SL_N},M)\to \cdots$$ is one of 
$A^{\SL_N}$-modules,
and $A^{\SL_N}$ is noetherian by invariant theory.
Also consider the long exact sequence 
$$\cdots\to H^i(\SL_N,k[\SL_N/U]\otimes M')\to H^i(\SL_N,k[\SL_N/U]\otimes
M)\to \cdots$$ \qed

More generally one has

\begin{Lemma}
Let 
$0\to M_0\to M_1\to \cdots\to M_q\to 0$ be a complex of \AG-modules
whose homology modules $\ker(M_i\to M_{i+1})/\im(M_{i-1}\to M_{i})$
are in $\N$, for $i=0,\dots,q$. If $q$ of the $M_i$ are
in $\N$, so
is the last one.
\end{Lemma}

\paragraph{Proof}
This is a routine consequence of the two out of three  property. \qed

\section{Picard graded Cox rings}\label{section 4}
If $V$ is a finite dimensional $k$-vector space, we denote its dual by $V^\#$.
For $1\leq s\le N$, let $\Gr(s)$ be the Grassmannian parametrizing 
$s$-dimensional subspaces of the 
dual $\nabla(\varpi_1)^\#$ of the
defining representation of $\GL_N$. Let $\cO(1)$ denote as usual
the ample generator of the Picard group of $\Gr(s)$. We wish to view it as
a $G$-equivariant sheaf. To this end consider the parabolic subgroup
$P=\{\;g\in G\mid g_{ij}=0\mbox{ for }i>N- s,~j\leq N- s\;\}$ and identify $\Gr(s)$
with $G/P$. Then a $G$-equivariant vector bundle is the associated bundle
of its fiber over $P/P$, where this fiber is a $P$-module. For the line bundle
$\cO(1)$ we let $P$ act by the weight $\varpi_N-\varpi_{N-s}$ on the fiber
over $P/P$. With this convention
$\Gamma(\Gr(s),\cO(1))$ is the Schur module $\nabla(\varpi_s)$, cf.\
\cite[II 2.16]{Jantzen}. 
More generally, for $n\geq0$ one has $\Gamma(\Gr(s),\cO(n))=\nabla(n\varpi_s)$.
So $$A\langle s\rangle =\bigoplus_{n\geq0}\Gamma(\Gr(s),\cO(n))$$
is a good $G\Z$-algebra. Recall that $\Delta=\Z^r$.
Let $1\leq s_i\le N$ be given for $1\leq i\leq r$.
Then the Cox ring $A\langle s_1\rangle \otimes\cdots\otimes A\langle s_r\rangle $ of 
$\Gr(s_1)\times\cdots\times \Gr(s_r)$ is a good \GD-algebra.
We put $C=C_0\otimes A\langle s_1\rangle \otimes\cdots\otimes A\langle s_r\rangle $ ,
where $C_0$ is a polynomial algebra on finitely many generators with trivial
$G$-action, and  $C_0$ is placed in degree zero. 
Then $C$ is also a good \GD-algebra. We wish to prove

\begin{Proposition}\label{Picardgraded}
Every \CGD-module is negligible.
\end{Proposition}

The proof will be by induction on the rank $r$ of $\Delta$.
It will be finished in \ref{endproofPicardgraded}.
As base of the induction we use

\begin{Lemma}
A \CG-module $M$ that is noetherian over $C_0$ is negligible. 
\end{Lemma}
\paragraph{Proof}
(Taken from \cite{cohGrosshans}.)
As  $M$ is a finitely generated $C_0$-module
it has only finitely many weights. Therefore the argument used
in \cite{Friedlander-Parshall} to show that finite dimensional $G$
modules have finite good filtration dimension, applies to $M$.

As ${\SL_N}$ is reductive, it is well known 
\cite[Thm. 16.9]{Grosshans book} that
$H^0({\SL_N},M)$ is a finite $C_0^{\SL_N}$-module. So we argue by dimension shift.
As $M$ has only finitely many weights,
one may choose $s$ so large that all weights of
$M\otimes k_{-(p^s-1)\rho}$ are anti-dominant, where $\rho=\sum_{i=1}^{N-1}\varpi_i$.
Let $\st_s$ denote the $s$-th Steinberg module $\ind_{B^+}^G(k_{-(p^s-1)\rho})$.
Then
$M\otimes \st_s=\ind_{B^+}^G(M\otimes k_{-(p^s-1)\rho})$ has by Kempf vanishing
a good filtration
and
therefore $M\otimes \st_s\otimes \st_s$ has a good filtration \cite[II 4.21]{Jantzen}.
Then $H^i({\SL_N},M)$ is the cokernel of $H^{i-1}({\SL_N},M\otimes \st_s\otimes \st_s)
\to H^{i-1}({\SL_N},M\otimes \st_s\otimes \st_s/M) $ for $i\geq1$.\qed

\begin{Notation}
For $1\leq q\leq r$ we denote by $C^{\widehat  q}$ the subring 
$\bigoplus_{I_q=0}C_I$.
\end{Notation}

We further assume $r\geq1$. The inductive hypothesis then gives:
\begin{Lemma}\label{inductive assumption}Let $1\leq q\leq r$.
If the \CGD-module $M$ is noetherian over the subring $C^{\widehat  q}$,
then $M$ is negligible.
\end{Lemma}

\section{Coherent sheaves}
We now have $\Proj(C)=\Spec(C_0)\times \Gr(s_1)\times\cdots\times \Gr(s_r)$.
Call the projections of $\Proj(C)$ onto its respective factors 
$\pi_0$, \ldots , $\pi_r$.
For $I\in\Delta$ define the coherent sheaf 
$\cO(I)=\bigotimes_{i=1}^r\pi_i^*(\cO(I_i))$.
So $C=\bigoplus_{I\geq0}\Gamma(\Proj(C),\cO(I))$.
For a \CGD-module $M$ let $M^\sim$ be the coherent
$G$-equivariant  sheaf \cite[2.1]{Brion Kumar}, cf.\  \cite[II F.5]{Jantzen},
on $\Proj(C)$
constructed as in \cite[II 5.1]{Hartshorne} from the $\Z$-graded
module $\diag(M):=\bigoplus_{I\in\diag(\Delta)}M_I$. 
Conversely, to a coherent sheaf $\cM$ on $\Proj(C)$, we associate
the $\Delta$-graded $C$ module 
$$\Gamma_*(\cM)=\bigoplus_{I\geq0}\Gamma(\Proj(C),\cM(I)),$$
where $\cM(I)=M\otimes\cO(I)$.
We also put $H_*^t(\cM)=\bigoplus_{I\geq0}H^t(\Proj(C),\cM(I))$.
Recall from \ref{gradings} that $E=(1,1,...,1)$,
 so that $\cO(E)$ is the natural very ample 
line bundle (relative to
$\Spec (C_0)$) on the Segre product of the Grassmannians in Pl\"ucker
embeddings.

\begin{Lemma}
If $\cM$ is a $G$-equivariant coherent sheaf  on $\Proj(C)$,
then the  $H_*^t(\cM)$ are \CGD-modules.
\end{Lemma}
\paragraph{Proof}
So we have to show that $H_*^t(\cM)$ is noetherian as a
$C$-module.
This is clear for $t>\dim(\Proj(C))$, so we argue by descending
induction on $t$. Assume the result for all larger values of $t$.
By Kempf vanishing 
$\bigoplus_{q\geq0}\bigoplus_{n\geq0}H^q(Gr(s),\cO(i+n))$
is a noetherian $\bigoplus_{n\geq0}\Gamma(Gr(s),\cO(n))$
module, for any $i\in\Z$,
so by a K\"unneth theorem 
$\bigoplus_{q\geq0}H_*^q(\Proj(C),\cO(I))$ is a noetherian $C$-module
for any $I\in \Delta$.
Now write $\cM$ as a quotient of some $\cO(iE)^a$ and use
the long exact sequence
$$\cdots\to H_*^t(\cO(iE)^a)\to H_*^t(\cM)\to H^{t+1}
(\dots)\to\cdots$$
to finish the induction step.\qed

\begin{Notation}
If $M$ is a $\Delta$-graded module and $I\in\Delta$, then $M(I)$ is the
$\Delta$-graded module with $M(I)_J=M_{I+J}$.
Further $M_{\geq I}$ denotes $\bigoplus_{J\geq I}M_J$.
\end{Notation}

\begin{Lemma}\label{diagtest}
If $I\geq0$, then the ideal $C_{\geq I}$ of $C$ is generated by $C_I$.\\
If $M$ is a \CGD-module with $M_{nE}=0$ 
for $n> >0$, then  $M_{\geq nE}=0$ for $n> >0$.
\end{Lemma}

\paragraph{Proof}
The ideal is generated by $C_I$ because $C$ is generated over $C_0$ by the
$C_{e_i}$.
Let $m\in M_I$. Choose $J\geq 0$ with $I+J\in\diag(\Delta)$. Then $mC_{J+qE}$
vanishes for $q\gg 0$, so $(mC)_{\geq I+J+qE}=0$ for $q\gg 0$.
Now use that $M$ is finitely generated over $C$.
\qed

\begin{Lemma}
If $M$ is a \CGD-module, then there is an $n_0$
so that if $I=nE=(n,\dots,n)\in\Delta$ with $n>n_0$,
then $M_{\geq I}=\Gamma_*(M^\sim)_{\geq I}$.
\end{Lemma}

\paragraph{Proof}
Recall \cite[II Ex.5.9]{Hartshorne} that we have a natural map
$\diag(M)\to \diag(\Gamma_*(M^\sim))$
whose kernel and cokernel live in finitely many degrees.
Consider the maps $f:\diag(M)\otimes_{\diag(C)}C\to M$
and $g:\diag(M)\otimes_{\diag(C)}C\to \Gamma_*(M^\sim)$.
If $N$ is the kernel or cokernel of $f$ or $g$ then $N_{nE}=0$ for 
$n\gg 0$. Now apply the previous lemma.
\qed

\begin{Lemma}
If $M$ is a \CGD-module and $I\in\Delta$, then $M/M_{\geq I}$ is
negligible.
\end{Lemma}
\paragraph{Proof}
As $M$ is finitely generated over $C$, there is $J<I$ with
$M=M_{\geq J}$. Now note that for $1\leq q\leq r$ and $K\in \Delta$
the module $M_{\geq K}/M_{\geq K+e_q}$ is negligible by \ref{inductive assumption}.
\qed

\begin{Definition}
In view of the above we call an equivariant coherent sheaf $\cM$
on $\Proj(C)$
negligible when $\Gamma_*(\cM)$ is negligible.
\end{Definition}

The following Lemma is now clear:

\begin{Lemma}Let $I\in\Delta$.
A $G$-equivariant coherent sheaf $\cM$ on $\Proj(C)$ is negligible if and only
if $\cM(I)$ is negligible.
\end{Lemma}

\begin{Lemma}
Let $$0\to\cM'\to\cM\to\cM''\to0$$
be an exact sequence of $G$-equivariant coherent sheaves on $\Proj(C)$.
There is $I\in\Delta$ with 
$$0\to\Gamma_*(\cM')_{\geq I}\to\Gamma_*(\cM)_{\geq I}
\to\Gamma_*(\cM'')_{\geq I}\to0$$
exact.
\end{Lemma}
\paragraph{Proof}The line bundle $\cO(E)$ is ample. 
Apply Lemma \ref{diagtest} to the homology sheaves of the
complex $$0\to\Gamma_*(\cM')\to\Gamma_*(\cM)
\to\Gamma_*(\cM'')\to0.$$
\qed

\begin{Lemma}\label{half}
For every $I\in\Delta$ the sheaf $\cO(I)$ is negligible.
If $\cF$ is a $G$-equivariant coherent sheaf on $\Proj(C)$ so that
$\Gamma_*(\cF)$ has finite good filtration dimension, then $\cF$ is negligible.
\end{Lemma}
\paragraph{Proof}
The first statement follows from the fact that $C$ is negligible.
As for the second, there is an equivariant exact sequence, 
$$0\to \cE\to \cO(i_qE)\otimes V_q\to\cdots\to\cO(i_1E)\otimes V_1\to\cF\to 0$$
with $\cE$ a vector bundle, and each $V_i$ a finite dimensional $G$-module. 
Note that $\Gamma_*(\cE)_{\geq nE}$ has finite good
filtration dimension for $n\gg 0$.  
Let $d=\lim_{n\to\infty}\dim_\nabla(\Gamma_*(\cE)_{\geq nE})$.
If $d=0$ then some $\Gamma_*(\cE)_{\geq J}$ has no higher ${\SL_N}$-cohomology and
is thus negligible by invariant theory \cite[Thm. 16.9]{Grosshans book}.
So we argue by induction on $d$. Say $d>0$.
As $\cE$ is a vector bundle, there is
short exact sequence of equivariant
vector bundles $0\to\cE\to \cO(nE)\otimes V\to \cE'\to 0$,
with $V$ a finite dimensional $G$-module. Any finite dimensional $G$-module
can be embedded into one with good filtration by \cite{Friedlander-Parshall},
so we may assume $V$ has good filtration.
As $\cE'$ has a smaller $d$ \cite[Lemma 2.1]{cohGrosshans}, induction applies.
\qed

\section{Resolution of the diagonal}
We write $X=\Proj(C)$, $Y=\Proj(C^{\hat r})$, $Z=\Gr(s)$, where $s=s_r$.
So $X=Y\times Z$.
We now recall the salient facts from \cite{twisted}, \cite{diagonal}
 about the 
functorial resolution of the diagonal in $Z\times Z$.
As $Z$ is the Grassmannian that parametrizes the $s$-dimensional subspaces
of 
$\nabla(\varpi_1)^\#$, we have the tautological exact sequence of $G$-equivariant
vector bundles
on $Z$:
$$0\to\cS\to \nabla(\varpi_1)^\#\otimes \cO_Z\to \cQ\to 0,$$
where $\cS$ has as fiber above a point the subspace $V$
that the point parametrizes, and $\cQ$
has as fiber above this same point the quotient 
$\nabla(\varpi_1)^\#/V$.
Let $\pi_1$, $\pi_2$ be the respective projections $Z\times Z\to Z$.
Then the composite of the natural maps 
$\pi_1^*(\cS)\to \nabla(\varpi_1)^\#\otimes \cO_{Z\times Z}$
and 
$\nabla(\varpi_1)^\#\otimes \cO_{Z\times Z}\to\pi_2^*(\cQ)$ defines a section
of the vector bundle $\cHom(\pi_1^*(\cS),\pi_2^*(\cQ))$ whose zero scheme 
is the diagonal $\diag(Z)$ in $Z\times Z$. Dually, we get an exact sequence 
$\cHom(\pi_2^*(\cQ),\pi_1^*(\cS))\to \cO_{Z\times Z}\to \cO_{\diag Z}\to 0$,
where $\cO_{\diag Z}$ is the quotient by the ideal sheaf defining the diagonal.
As the rank $d$ of the vector bundle $\cE=\cHom(\pi_2^*(\cQ),\pi_1^*(\cS))$ equals
the codimension of $\diag(Z)$ in $Z\times Z$, the Koszul complex
$$0\to \bigwedge^d\cE\to\cdots\to\cE\to\cO_{Z\times Z}\to \cO_{\diag Z}\to 0$$
is exact. Now each $\bigwedge^i\cE$ has a finite filtration whose associated
graded is 
$$\bigoplus S^\alpha \pi_1^*(\cS)\otimes (S^{\tilde\alpha}\pi_2^*(\cQ))^\#,$$
where $\alpha $ runs over partitions of $i$ with at most $\rank(\cS)$ parts,
so that moreover the conjugate partition $\tilde \alpha$ has at most $\rank(\cQ)$
parts.

\paragraph{Plan}Now the plan is this: Let $\pi_{1,2}$ be the projection
of  $Y\times Z\times Z$  onto the product $Y\times Z$ of the first two factors,
let $\pi_2$ be the projection onto the middle factor $Z$, and so on.
If $M$ is a \CGD-module, tensor the pull-back along $\pi_{2,3}$ of the
 Koszul complex with
$\pi_{1,3}^*(M^\sim)$, take a high Serre twist
and then the direct image along $\pi_{1,2}$ to $X$. 
On the one hand $(\pi_{1,2})_*(\pi_{1,3}^*(M^\sim)\otimes \cO_{\diag Z})$
is just $M^\sim$, but on the other hand
the salient facts above allow us to express it in
terms of negligible  \CGD-modules. This will prove that $M$ is negligible.
We now proceed with the details.

\begin{Remark}
Instead of functorially resolving the diagonal in $Z\times Z$, we could have
functorially resolved the diagonal in $X\times X$.
\end{Remark}

\begin{Notation}
On a product like $Y\times Z$ an exterior tensor product
$\pi_1^*(\cF)\otimes\pi_2^*(\cM)$ is denoted $\cF\boxtimes\cM$.
\end{Notation}

\begin{Lemma}
Let $\cF$ be a $G$-equivariant coherent sheaf on $Y$, and $\alpha$ a partition of
$i$ with at most $s$ parts, $i\geq0$.
The sheaf $\cF\boxtimes S^\alpha(\cS)$ on $X=Y\times Z$ is
negligible.
\end{Lemma}
\paragraph{Proof}
By the inductive assumption 
$$\Gamma_*(\cF)=\bigoplus_{I\in\Z^{r-1},~I\geq0}
\Gamma(Y,\cF(I))$$
is a
$C^{\widehat r}$-module
with finite good filtration dimension.
The vector bundle  $\cS$ on $Z=G/P$ is associated with the
irreducible $P$-representation with lowest weight $-\epsilon_{N-s+1}$.
This representation may be viewed as $\ind_{B^+}^P(-\epsilon_{N-s+1})$,
where $-\epsilon_{N-s+1}$ also stands for
the one dimensional $B^+$ representation with weight $-\epsilon_{N-s+1}$.
Say $\rho:P\to P^-$ is the isomorphism which sends
a matrix to its transpose inverse. Then 
$\ind_{B^+}^P(-\epsilon_{N-s+1})=
\rho^*\ind_{B^-}^{P^-}(\epsilon_{N-s+1})$.
One finds that $S^\alpha(\cS)$ is associated with
$\rho^*\ind_{B^-}^{P^-}\left(\sum_i\alpha_i\epsilon_{N-s+i}\right)=
\ind_{B^+}^P\left(-\sum_i\alpha_i\epsilon_{N-s+i}\right)$.
(This is the rule $S^\alpha(\nabla_{\GL_s}(\varpi_1))=\nabla_{\GL_s}(\alpha)$
in disguise.)
Then $S^\alpha(\cS)(n)$ is associated with 
$\ind_{B^+}^P\left(-\sum_i\alpha_i\epsilon_{N-s+i}+n\varpi_N-n\varpi_{N-s}\right)$.
For $n\geq \alpha_1$ the weight 
$-\sum_i\alpha_i\epsilon_{N-s+i}+n\varpi_N-n\varpi_{N-s}$ is an anti-dominant
polynomial weight, so 
$\sum_{n\geq\alpha_1}\Gamma(Z,S^\alpha(\cS)(n))$
has a good filtration by transitivity of induction \cite[I 3.5, 5.12]{Jantzen}.
Then $\Gamma_*(\cF\boxtimes S^\alpha(\cS))_{\geq I}$
has finite good filtration dimension \cite[Lemma 2.1]{cohGrosshans}
for $I=(0,\dots,0,\alpha_1)$
and the result follows from Lemma \ref{half}.\qed

\begin{Assumption}Recall we are trying to prove that $M$ is negligible.
As in the proof of Lemma \ref{half}, we may reduce to the case that $M^\sim$ is a
vector bundle. We further assume this.
\end{Assumption}

\begin{Lemma}
For $n\gg 0$ the sheaf 
$$(\pi_{12})_*\left(\pi_{13}^*(M^\sim)\otimes\bigg(\cO(nE)\boxtimes\cO(n)\bigg)
\otimes\pi_{23}^*(\bigwedge^i\cE)\right)$$ is
negligible.
\end{Lemma}
\paragraph{Proof}The sheaf $\cO(E)\boxtimes\cO(1) $ is ample.
So \cite[Thm. 8.8]{Hartshorne} the sheaf in the Lemma
has a filtration with layers of the form
$$(\pi_{12})_*\left(\pi_{13}^*(M^\sim)\otimes\bigg(\cO(nE)\boxtimes\cO(n)\bigg)
\otimes\pi_{23}^*\bigg( S^\alpha (\cS)\boxtimes \cG\bigg)\right).$$
Say $f:Y\times Z\to Y$ is the projection.
Now use $(\pi_{12})_*\circ\pi_{13}^*=f^*\circ f_*$ and
a projection formula for $(\pi_{12})_*$ to rewrite the layer
in the form $(\cF\boxtimes S^\alpha(\cS))(I)$
for some $I\in\Delta$, with $I$ depending on $n$.\qed

\paragraph{End of proof of Proposition \ref{Picardgraded}}
Proposition \ref{Picardgraded} now follows from

\begin{Lemma}\label{endproofPicardgraded}
$M^\sim$ is negligible.
\end{Lemma}
\paragraph{Proof}
{}From the Koszul complex and the previous Lemma we conclude \cite[Thm. 8.8]{Hartshorne}
that for $n\gg 0$ 
the sheaf $$(\pi_{12})_*\left(\pi_{13}^*(M^\sim)\otimes\bigg(\cO(nE)\boxtimes\cO(n)\bigg)
\otimes\pi_{23}^*(\cO_{\diag(Z)})\right)$$ is
negligible. This sheaf equals $M^\sim(I)$ for some $I\in \Delta$.
\qed

\section{Differently graded Cox rings}
Let $c:\{1,\dots,r\}\to\{1,\dots,q\}$ be surjective.
Put $\Lambda=\Z^q$. We have a \emph{contraction map}, also denoted $c$,
from $\Delta$ to $\Lambda$ with $c(I)_j=\sum_{i\in c^{-1}(j)}I_i$.
Through this contraction we can view our $\Delta$-graded $C$ as $\Lambda$-graded.
We now have the following generalization of Proposition \ref{Picardgraded}:

\begin{Proposition}\label{retractedgraded}
Every \CGL-module is negligible.
\end{Proposition}
This will be proved by descending induction on $q$, with fixed $r$.
The case $q=r$ is clear. So let  $q<r$ and assume
the result for larger values of $q$.
We may assume $c(r-1)=c(r)=q$. (Otherwise rearrange the factors.)
Recall $X=\Proj(C)$, $X=Y\times Z$, with $Y=\Proj(C^{\widehat r})$, 
$Z=\Proj(A\langle s\rangle)$.

\begin{Notation}
Let $\m$ be the irrelevant maximal ideal $\bigoplus_{i>0}A\langle s\rangle _i$
of $A\langle s\rangle$. If $M$ is a \CGL-module, put $M_{\geq i}=\m^iM$, and
$\gr^iM=M_{\geq i}/M_{\geq i+1}$. If $I\in \Lambda$, put $(M_I)_{\geq i}=M_I\cap \m^iM$,
and $\gr^iM_I=(M_I)_{\geq i}/(M_I)_{\geq i+1}$.
We put a $\Z^{q+1}$-grading on $\gr M=\bigoplus_i\gr^i M$ with 
$$(\gr M)_I=\gr^{I_{q+1}}M_{(I_1,\dots,I_{q-1},I_q+I_{q+1})}.$$
In particular all this applies when $M=C$. Then $\gr C$ may be identified with $C$
and the $\Z^{q+1}$-grading on $\gr C$ is a contracted grading to which the inductive
assumption applies. Write $\Phi=\Z^{q+1}$. Then $\gr M$ is a $CG\Phi$-module.
\end{Notation}

Let $M$ be a \CGL-module.
By the inductive assumption $\gr M$ has finite good filtration dimension and
each $H^i({\SL_N},\gr M)$ is a noetherian $(\gr C)^{\SL_N}$-module.
We still have to get rid of the grading.
The filtration $M_{\geq 0}\supseteq M_{\geq 1}\cdots$ induces a filtration
of the Hochschild complex \cite[I.4.14]{Jantzen} whence a spectral sequence
$$E(M):E_1^{ij}=H^{i+j}({\SL_N},\gr^{i}M)\Rightarrow H^{i+j}({\SL_N}, M).$$
It lives in two quadrants. The spectral sequence $E(M)$ is a direct sum
of spectral sequences $E(M_I)$, $I\in\Lambda$. As each $M_I$ has a finite
filtration, each $E(M_I)$ stops, meaning that there is an $a$ so that the
differentials in $E_b^{**}(M_I)$ vanish for $b\geq a$. Thus 
$E_a^{**}(M_I)=E_\infty^{**}(M_I)$ is an associated graded of the abutment $H^*({\SL_N},M_I)$.

\begin{Lemma}
$E(M)$ also stops and its abutment is a noetherian $C^{\SL_N}$-module.
\end{Lemma}
\paragraph{Proof}
The spectral sequence
$E(C)$
is pleasantly boring: It does not just degenerate,
even its abutment is the same as its $E_1$.
The spectral sequence $E(M)$ is a module over it
\cite[Theorem 3.9.3]{Benson II}, \cite{Massey}.
In particular, $E(M)$ is a module over $C^{\SL_N}$. But $E_1^{**}(M)$ is noetherian
over $C^{\SL_N}=(\gr C)^{\SL_N}$.
So the usual argument (see \cite[Lemma 3.9]{reductive} or
\cite[Lemma 7.4.4]{Evens book})
shows that $E(M)$ stops and that $E_\infty^{**}(M)$ is noetherian over $C^{\SL_N}$.
As the filtrations on the abutments of the $E(M_I)$ are finite, it follows
that the abutment of $E(M)$ is finitely generated over $C^{\SL_N}$.
\qed

\begin{Lemma}
$M$ has finite good filtration dimension.
\end{Lemma}
\paragraph{Proof}
As each $M_I$ is finitely filtered, $\dim_\nabla(M_I)\leq\dim_\nabla(\gr M_I)$.
\qed\\

This finishes the proof of Proposition \ref{retractedgraded}.

\section{Variations on the Grosshans grading}\label{G=SL}
In this section we will be concerned with representations of $\SL_N$.
Mutatis mutandis everything also
applies to other connected reductive groups. We now write $G=\SL_N$,
with subgroups $B^+$, $B^-$, $T$, $U$ defined in the usual manner. (So they are now
the intersections with $\SL_N$ of the subgroups of $\GL_N$ that had these names.)
As explained above, the Grosshans graded $\gr V$ of an $\SL_N$-module $V$ has
a $\Z$-grading. We also need a $\Lambda$-graded version, where $\Lambda$ is 
the weight lattice of $\SL_N$. In \cite{vdkallen book} such a
version was studied using a total order on weights known as the length-height order.
It was claimed incorrectly in \cite{cohGrosshans} that one might as well use the dominance
order which is only a partial order. And it was claimed incorrectly in \cite{cohGrosshans}
that the resulting $\SL_N$-module is isomorphic with $\gr V$.
Both claims are correct when $V$ has good filtration, but they
are wrong in general. See example \ref{not so} below.
The claims are repeated in \cite{reductive}, \cite{fgfdim}.
Let us now introduce a $\Lambda$-graded version that is closer to the Grosshans
graded than the version based on length-height order. (Length-height order was 
appropriate
when dealing with the category of $\SL_N$-modules as embedded into the
larger category of $B$-modules.)
Following Mathieu \cite{Mathieu G}
we choose a second linear height function $E:\Lambda\otimes \R\to \R$
with $E(\alpha)>0$ for every positive root $\alpha$, but now with $E$ injective
on $\Lambda$. We
define a total order on weights by first ordering them by Grosshans height,
then for fixed Grosshans height by $E$. With this total order, denoted $\leq$,  we put:

\begin{Definition}
If $V$ is a $G$-module, and $\lambda$ is  a weight,
then $V_{\leq\lambda}$ denotes the largest $G$-submodule all whose weights
$\mu$ satisfy $\mu\leq\lambda$ in the total order.
For instance, $V_{\leq0}$ is the module of invariants
$V^G$. Similarly $V_{<\lambda}$ denotes the largest $G$-submodule all 
whose weights
$\mu$ satisfy $\mu<\lambda$.  Note that $V\mapsto V_{\leq \lambda}$ is a truncation functor for a 
saturated set
of dominant weights \cite[Appendix A]{Jantzen}.
 So this functor fits in the usual
highest weight category picture.
As in \cite{vdkallen book}, we form the
$\Lambda$-graded module
$$\gr_{\Lambda} V=\bigoplus_{\lambda\in \Lambda}V_{\leq\lambda}/V_{<\lambda}.$$
Each $\gr_\lambda V=V_{\leq\lambda}/V_{<\lambda}$
 has a $B^+$-socle
$(\gr_\lambda V)^U=V^U_\lambda$ of weight
$\lambda$. We always view $V^U$ as a $B^-$-module through restriction
(inflation)
along the
homomorphism $B^-\to T$.
Then $\gr_\lambda V$ embeds naturally in its `good filtration
hull' $\hull_\nabla(\gr_\lambda V)
=\ind_{B^-}^GV^U_\lambda$.
This good filtration hull has the same $B^+$-socle.
\end{Definition}

If $\lambda$ is not dominant, then $\gr_\lambda V$ vanishes, because
its socle vanishes.
Note that $\bigoplus_{\hgt(\lambda)=i}\gr_{\lambda} V$ is the associated graded of
a filtration of $\gr_i V$, where $\gr_{\lambda} V$ refers to a graded component of $\gr_{\Lambda} V$
and $\gr_i V$ to one of $\gr V$. Both $\gr_{\Lambda} V$ and $\gr V$ embed into the
good filtration hull $\ind_{B^-}^GV^U$, which is $\Lambda$-graded. But while 
$\gr_{\Lambda} V$ is a $\Lambda$-graded submodule of the hull, 
$\gr V$ need only be a $\Z$-graded submodule. Both $\gr_{\Lambda} V$ and $\gr V$
contain the socle of the hull.

\begin{Example}\label{not so} Take $p=2$, $N=3$. As group we may take $\SL_3$
or $\GL_3$.
Inside
$\nabla(3\varpi_1+\varpi_3)\oplus\nabla(3\varpi_2)$ take an indecomposable submodule $V$
of codimension one. Then $V$ has three composition factors.
It has a one dimensional head and its socle is the direct sum
of two irreducibles, whose highest weights have identical Grosshans height.
 It is easy to see that $\gr_{\Lambda} V$ has two indecomposable summands and
$\gr V$ just one. And using the dominance order as suggested in \cite{cohGrosshans}
would not even lead to an associated graded of $V$. The head gets lost.
\end{Example}

Although $\gr_\Lambda V$ need not coincide with $\gr V$ it shares some properties:

\begin{Lemma}\label{Lambda Grosshans}
\begin{enumerate}
\item If $A$ is a finitely generated $k$-algebra, so is $\gr_\Lambda A$.
\item If $A$ has good filtration, then $\gr_\Lambda A$ is isomorphic to $\gr A$ as
$k$-algebra.
\end{enumerate}
\end{Lemma}

\paragraph{Proof}Both $\gr A$ and $\gr_\Lambda A$ embed into 
  their good filtration hull $\ind_{B^-}^GA^U$, notation $\hull_\nabla(\gr A)$,
  cf. \cite[2.2]{cohGrosshans}. 
 The argument of Mathieu (see proof of \cite[Lemma 2.3]{cohGrosshans})
that this  $\hull_\nabla(\gr A)$ is the $p$-root closure
of $\gr A$ applies just as well to the subalgebra $\gr_\Lambda A$. 
Indeed it would even apply to the subalgebra $S$
of $\hull_\nabla(\gr A)$ generated by the socle of the hull. We argue as in the proof of
\cite[Theorem 9]{Grosshans contr}.
The finitely generated algebra 
$\hull_\nabla(\gr A)$ is integral over its finitely generated subalgebra $S$ and
$\gr_\Lambda A$ is an $S$-submodule of the hull. 
Then $\gr_\Lambda A$ must
be finitely generated. 
When $A$ has good filtration, $\gr_i A$ is already a direct sum of costandard
modules. So then passing to the associated graded of the filtration of $\gr_i A$ makes 
no difference. And the algebra structure on both $\gr A$ and $\gr_\Lambda A$
agrees with the algebra structure on the hull by \cite[Lemma 2.3]{fgfdim}.
\qed

\section{Proofs of the main results}
Let us now turn to the proof of Theorem \ref{maingood} for $\SL_N$.
Return to the notations introduced in section \ref{G=GL}. 
Thus $G=\GL_N$, with $T$ its maximal torus.
We assume the $\SL_N$-algebra $A$ has a good filtration and $M$ is a
noetherian $A$-module on which $\SL_N$ acts compatibly.
Put $\Lambda=\Z^{N-1}$ and identify $\Lambda$ with a sublattice of $X(T)$
by sending $\lambda\in\Lambda$ to $\sum_i\lambda_i\varpi_i$.
Also identify $\Lambda$ with $X(T\cap\SL_N)$ through the restriction
$X(T)\to  X(T\cap\SL_N)$. Thus a dominant $\lambda\in\Lambda$ gets identified with
a polynomial dominant weight.
For such $\lambda$
we may embed $\gr_\lambda A$ or
$\gr_\lambda M$ into its good filtration hull which is a direct
sum of restrictions to $\SL_N$ of the Schur module $\nabla_G(\lambda)$.
On the Schur module $\nabla_G(\lambda)$ the center of $G$ acts through
$\lambda$. This makes it natural to
use the $\Lambda$-grading on $\gr_\Lambda A$ and $\gr_\Lambda M$ 
to extend the action from $\SL_N$ 
to $\GL_N$, making the center of $\GL_N$ act through $\lambda$
on the graded pieces $\gr_\lambda A$ and $\gr_\lambda M$. We do that.
Next we imitate subsection 2.2 of \cite{fgfdim}.

\begin{Lemma}\label{extend}
Recall $A$ has a good filtration, so that $\gr_\Lambda A=\hull_\nabla(\gr_\Lambda A)$.
Let $R=\oplus_\lambda R_\lambda$ be a $\Lambda$-graded algebra with $G$-action
such that $R_\lambda=(R_\lambda)_{\leq \lambda}$. Then every $T$-equivariant
graded algebra homomorphism $R^U\to (\gr_\Lambda A)^U$ extends uniquely to a
$G$-equivariant graded algebra homomorphism $R\to \gr_\Lambda A$.
\end{Lemma}

\paragraph{Proof}
Use that $\hull_\nabla(\gr_\Lambda A)$  is an induced module.\qed\\

As the algebra
$(\gr_\Lambda A)^U=(\gr A)^U$ is finitely generated by Grosshans \cite{Grosshans contr}, 
it is also generated by finitely many
weight vectors. Consider one such weight vector $v$, say of weight $\lambda$.
Clearly $\lambda$ is dominant.
If $\lambda=0$, map a polynomial ring $P_v:=k[x]$ with trivial $G$-action to
 $\gr A$ by substituting $v$ for $x$. Also put $D_v:=1$.
Next assume $\lambda\neq0$. Let $\ell=N-1$ be the rank of $\Lambda$.
Recall the Cox rings $A\langle i\rangle$ of section \ref{section 4}. Define a $T$-action
on the $\Lambda$-graded algebra
$$P=\bigotimes_{i=1}^{\ell}A\langle i\rangle$$
by letting $T$ act on $\bigotimes_{i=1}^{\ell}\Gamma(\Gr(i),\cO(m_i))$
through weight $\sum_im_i\varpi_i$. So now we have a $G\times T$-action on $P$,
and the $T$-action corresponds with the $\Lambda$-grading.
Observe
that by  the tensor product
property \cite[Ch.~G]{Jantzen}
the algebra $P$ has a good filtration for the $G$-action.
Let $D$ be the scheme
theoretic kernel of $ \lambda$.
So $D$ has character group
$X(D)=X(T)/\Z \lambda$ and $D=\Diag(X(T)/\Z \lambda)$ in the notations of
\cite[I.2.5]{Jantzen}.
The subalgebra $P^{1\times D}$ is a graded algebra with
good filtration such that its subalgebra
$P^{U\times D}$ contains a polynomial algebra on one generator $x$ of weight
$\lambda\times \lambda$. In fact, this polynomial subalgebra contains all the
weight vectors in $P^{U\times D}$ whose weight is of the form $\nu\times \nu$.
 The other weight vectors in $P^{U\times D}$
have weight of the form $\mu\times\nu$
with $\nu$ an integer multiple of $\lambda$ and $\mu<\nu$. These other weight vectors
span an ideal in $P^{U\times D}$.
By lemma \ref{extend} one easily constructs
 a $G$-equivariant algebra homomorphism
$P^{1\times D}\to \gr_\Lambda A$ that maps $x$ to $v$.
Write it as $P_v^{1\times D_v}\to \gr_\Lambda A$, to stress the dependence on $v$.

The direct product $D$ of the $D_v$ is a diagonalizable group. It acts on
the tensor product $C$ of the finitely many $P_v$. This $C$ is $\Lambda$-graded.
We have a  graded algebra map $C^D\to \gr_\Lambda A$.
It is surjective because its image has good filtration
(\cite[Ch.~A]{Jantzen}) and contains $(\gr A)^U$.
We have proved

\begin{Lemma}There is a graded $G$-equivariant surjection $C^D\to \gr_\Lambda A$,
where the
$G\times D$-algebra $C$  is a good $G\Lambda$ algebra as in \ref{retractedgraded}. 
\end{Lemma}

Now recall $M$ is a noetherian $A$-module on which $G$ acts compatibly, meaning
that the structure map $A\otimes M\to M$ is a map of $G$-modules.
Form the `semi-direct product ring' $A\ltimes M$ whose underlying $G$-module
is $A\oplus M$, with product given by
$(a_1,m_1)(a_2,m_2)=(a_1a_2,a_1m_2+a_2m_1)$. By \ref{Lambda Grosshans} $\gr_\Lambda(A\ltimes M)$ is
a finitely generated algebra, so we get

\begin{Lemma}
$\gr_\Lambda M$ is a noetherian $\gr_\Lambda A$-module.
\end{Lemma}

This is of course very reminiscent of the proof of the lemma
\cite[Thm. 16.9]{Grosshans book} telling that
$M^G$ is a noetherian module over the finitely generated $k$-algebra
$A^G$. We  will tacitly use
its counterpart for diagonalizable actions, cf. \cite{Borsari-Santos},
\cite[I.2.11]{Jantzen}.

Now this lemma implies that $C\otimes_{C^D}\gr_\Lambda M$ is a \CGL-module, so by Proposition 
\ref{retractedgraded} the following analogue of \cite[Lemma 2.7]{fgfdim} holds.

\begin{Lemma}
$C\otimes_{C^D}\gr_\Lambda M$ has finite good filtration dimension and each 
$H^i(\SL_N,C\otimes_{C^D}\gr_\Lambda M)$ is a noetherian $C^{\SL_N}$-module.
\end{Lemma}

\begin{Remark}Note that $C\otimes_{C^D}\gr_\Lambda M$ actually has finite
Schur filtration dimension. Indeed we only need Proposition 
\ref{retractedgraded} for polynomial \CGL-modules. On the other hand
the reader may prefer to prove a version of Proposition 
\ref{retractedgraded} for $\SL_N$ rather than extending the action 
on $\gr_\Lambda A$ and $\gr_\Lambda M$ from $\SL_N$ to $G=\GL_N$. We now
have to restrict back to $\SL_N$ anyway.
\end{Remark}

Now we get  the analogue of \cite[Lemma 2.8]{fgfdim}
 
\begin{Lemma}
The module $\gr_\Lambda M$ has finite good filtration dimension and 
$\oplus_iH^i(\SL_N,\gr_\Lambda M)$ is a noetherian $A^{\SL_N}$-module.
\end{Lemma}
\paragraph{Proof}
Extend the $D$-action on $C$ to $C\otimes_{C^D}\gr_\Lambda M$ by using the trivial
action on the second factor. Then we have a $G\times D$-module structure
on $C\otimes_{C^D}\gr_\Lambda M$. As $D$ is diagonalizable,
$C^D$ is a direct summand of $C$ as a $C^D$-module \cite[I.2.11]{Jantzen}
and $(C\otimes_{C^D}\gr_\Lambda
M)^{1\times D}=\gr_\Lambda M$ is a direct summand of the $G$-module
$C\otimes_{C^D}\gr_\Lambda M$.
It follows that $\gr_\Lambda M$ also has finite good filtration dimension and it
follows that each $H^i(\SL_N,C\otimes_{C^D}\gr_\Lambda M)^{1\times D}=H^i(\SL_N,\gr_\Lambda M)$
is a noetherian $C^{\SL_N\times D}$-module. And there are only finitely many $i$
for which $H^i(\SL_N,\gr_\Lambda M)$ is nonzero.
But the action of $C^{\SL_N\times D}$ on $\gr_\Lambda M$
factors through $(\gr_\Lambda A)^{\SL_N}$, so we see
that each $H^i(\SL_N,\gr_\Lambda M)$ is a noetherian $(\gr_\Lambda A)^{\SL_N}$-module.
And one always has $(\gr_\Lambda A)^{\SL_N}=(\gr_0 A)^{\SL_N}=A^{\SL_N}$.
\qed

\paragraph{End of proof of Theorems \ref{maingood}, \ref{mainSchur}.}
We see that each $\gr_\lambda M$
is negligible as $(A^{\SL_N})G$-module. 
Enumerate the dominant weigths in
$\Lambda$ as
$\lambda_0$,  $\lambda_1$, \dots\ according to our total order on weights.
Note there are only finitely many dominant weights of given Grosshans height in $\Lambda$,
so that the order type of the set of dominant weights in $\Lambda$
is indeed just that of $\mathbb N$. (This would be false for the set of
dominant weights in $X(T)$.)
Using the two out of three property \ref{two/three} we see by induction that
$M_{\leq\lambda_n}$ is negligible as $(A^{\SL_N})G$-module. 
Moreover, as 
$\bigoplus_{i, \mu}H^i(\SL_N,\gr_\mu M)$ is noetherian over $A^{\SL_N}$,
there are only finitely many nonzero
$H^i(\SL_N,\gr_\mu M)$.
 So by a limit argument \cite[I Lemma 4.17]{Jantzen}
 each $H^i(\SL_N,M)$ is 
 a noetherian $A^{\SL_N}$-module.
 There is an $m$ with
$H^m(\SL_N,k[\SL_N/U]\otimes\gr_\lambda M)=0$ for all $\lambda\in\Lambda$. 
So by a similar limit argument $H^m(\SL_N,k[\SL_N/U]\otimes M)=0$ and
 $M$ has finite good filtration dimension.
This proves the theorem for the $\SL_N$ case. The $\GL_N$ case follows from
the $\SL_N$ case, using that $H^i(\GL_N,M)=H^i(\SL_N,M)^\Gm$ for a $\GL_N$-module
$M$.
Of course Theorem \ref{mainSchur} follows from Theorem \ref{maingood}
by Proposition \ref{Schurgood}.

\paragraph{Proof of Corollary \ref{mainGrosshans}}
Now let $A$ be any finitely generated commutative $k$-algebra on which $\SL_N$
acts rationally by $k$-algebra automorphisms. We argue as in the proof of
\cite[Proposition 3.8]{cohGrosshans}. Recall again the following result of Mathieu \cite{Mathieu G},
cf.~\cite[Lemma 2.3]{cohGrosshans}
\begin{Lemma}
For every $x\in\hull_\nabla(\gr A)$, there is an integer
$r\geq0$, so that $x^{p^r}\in \gr A$.
\end{Lemma}
But $\hull_\nabla(\gr A)$ is finitely generated by Grosshans, 
so let us fix $r$ so that for every $x\in\hull_\nabla(\gr A)$, one has $x^{p^r}\in \gr A$.
By \cite[Theorem 1.5, Remark 1.5.1]{Friedlander-Suslin}
the
ring $R=H^*(G_r,\gr A)^{(-r)}$ is a finite module
over the algebra $$\bigotimes_{a=1}^rS^*((\gl_n)^\#(2p^{a-1}))\otimes
\hull_\nabla(\gr A).$$ This algebra has a good filtration by \cite[4.3]{Andersen-Jantzen}, 
\cite[Chapter G]{Jantzen}. By Theorem \ref{maingood}
the ring $R$ has finite good filtration dimension. 
Therefore there are only finitely many $i$ with $E_2^{i*}\neq 0$
in the spectral sequence
$$E_2^{ij}=H^i(G/G_r,H^j(G_r,\gr A))\Rightarrow H^{i+j}(G,\gr A).$$ 
So this spectral sequence stops,
i.e. $E_s^{**}=E_\infty^{**}$ for some $s<\infty$.
By the same Theorem $H^*(G,R)$
is finite over the ring
$H^0(G,\bigotimes_{a=1}^rS^*((\gl_n)^\#(2p^{a-1}))\otimes
\hull_\nabla(\gr A))$, which is finitely generated by invariant theory \cite[Thm. 16.9]{Grosshans book}. 
So $H^*(G,R)=E_2^{**}$ is a finitely generated $k$-algebra.
Every page $E_a^{**}$ is a differential graded algebra in characteristic $p$,
so the $p$-th power of an even element passes to the next page.
Using this one sees that all pages are finitely generated as  $k$-algebras.
In particular, $E_\infty^{**}$ is finitely generated.
As the spectral sequence lives in the first quadrant, the abutment is also
finitely generated. \qed

\begin{Remark}
Similarly the $k$-algebra $H^*(\SL_N,\gr_\Lambda A)$ is finitely generated.
But $\gr_\Lambda A$ is even more graded than $\gr A$, and thus lies in the opposite 
direction of where we would like to go.
\end{Remark}
 
\phantomsection   
\end{document}